\documentclass[12pt]{amsart}

\usepackage{amssymb}
\usepackage{amsmath}
\usepackage{amscd}
\usepackage{float}
\usepackage{graphicx}
\usepackage{amsfonts}
\usepackage{pb-diagram}
\usepackage{eufrak}

\begin{document}

 \title{ An Affine Index  Polynomial Invariant of Virtual Knots}

\author{Louis H. Kauffman}
\address{ Department of Mathematics, Statistics and
 Computer Science, University of Illinois at Chicago,
 851 South Morgan St., Chicago IL 60607-7045, U.S.A.}
\email{kauffman@math.uic.edu}
\urladdr{http://www.math.uic.edu/$\tilde{~}$kauffman/}

\thanks{The author was partially supported by UIC}

\keywords{virtual knot, virtual link, writhe, Vassiliev invariant, quandle.}

\subjclass[2000]{57M27 }

\date{}
\maketitle


\begin{abstract}
This paper describes a polynomial invariant of virtual knots that is defined in terms of an integer labeling of the virtual knot diagram. This labeling is seen to derive from an essentially unique structure
of affine flat biquandle for flat virtual diagrams. The invariant is discussed in detail with many examples,including  its relation to previous invariants of this type and we show how to construct Vassiliev invariants from the same data.
 \end{abstract}

\section{Introduction}
This paper generalizes invariants of virtual knots defined by Z. Cheng \cite{Cheng} and by A. Henrich \cite{HenrichThesis}  to a new polynomial invariant of virtual knots and links. The invariant discussed herein is also related to the generalized parity invariants of H. Dye \cite{Dye}. In all these cases, an invariant is constructed in terms of {\it weights}, $W_{K}(c)$, associated to the crossings $c$ of an oriented  virtual knot $K$ and the invariants take the form of a polynomial defined by the equation
$$P_{K}(t) = \sum_{c} sgn(c)(t^{W_{K}(c)} -1)$$ where $sgn(c)$ denotes the sign of the crossing $c$ in  the oriented knot $K.$ The weights for the authors mentioned above are derived from the 
combinatorics of chord intersections in Gauss diagrams for the knots. In Cheng's case the polynomial 
utilizes parity and the weights are restricted to odd crossings (crossings corresonding to chords that intersect an odd number of other chords in the Gauss diagram). In Henrich's case the weights utilized are absolute values of the Cheng weights. The invariants in this paper are quite distinct from the index polynomial invariants discussed in \cite{NKamada,YH1,YH2}, while these invariants follow a similar pattern in the form of the polynomial. We call the polynomial invariant given in this paper the {\em Affine Index Polynomial} because in our approach the polynomial is a way of assembling a set of crossing
weights that are derived from a very simple affine biquandle structure on the underlying flat diagram.
This affine structure is explained in detail in the body of the paper.
\bigbreak

In this paper, we give a definition of weights $W_{K}(c)$ that exhibits them as differences related to an integer labeling of a flat virtual diagram associated with the knot diagram. This makes the system of weights very easy to compute and one can then develop the invariant on this basis. We call the resulting polynomial the {\it Affine Index Polynomial} $P_{K}(t)$ for an oriented virtual knot $K.$ The paper is organized as follows. In section 2 we review the definitions of virtual knot theory and flat virtual knot theory. In section 3 we give the definition of the Affine Index Polynomial and one example. In section 4 we prove the invariance of the polynomial. In section 5 we give a number of examples of computations of the invariant. We give examples where virtual knots with unit Jones polynomial are both shown to be non-classical and examples where this Index invariant cannot show that the knot is non-classical.
We give examples where invertibility is detected and we give examples where the non-triviality of a flat virtual knot is detected. We also give examples of cobordant knots and links that can be assigned an
invariant. In section 6 we show that the following formulas give Vasiliev invariants of order 
$\lceil n/2 \rceil$:
$$v_{n}(K) = (1/n!)(\sum_{c} sgn(c) W_{K}(c)^{n})$$  where $c$ runs over all the crossings of the knot 
$K$, $sgn(c)$ is the sign of the crossing, and $W_{K}(c)$ is the index of the crossing.
In section 8 we analyse linear affine flat biquandles and show that for a full biquandle structure, the 
integer labeling that we have used to define the invariant is essentially unique. It is a consequence of this analysis that there are other affine flat pre-biquandles (made explicit in this section) that can be used to make invariants of flat virtuals without the third Reidemeister move (virtual doodles). This will be the subject of a separate investigation.
\bigbreak

\noindent {\bf Acknowledgement.} It gives the author great pleasure to thank Lena Folwaczny for many conversations and to thank the Newton Institute for Mathematical
Sciences in Cambridge, UK for its hospitality in the final stages of the preparation of this paper.
\bigbreak

 \section{Recollection of Virtual Knot Theory}
\begin{figure}
     \begin{center}
     \includegraphics[width=8cm]{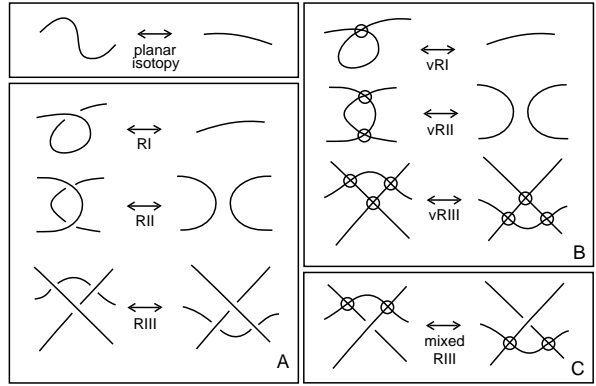}
     \caption{Virtual Isotopy}
     \label{virtualmoves}
\end{center}
\end{figure}

\begin{figure}
     \begin{center}
     \includegraphics[width=8cm]{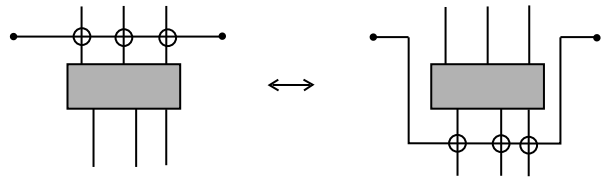}
     \caption{Detour Move}
     \label{detourmove}
\end{center}
\end{figure}

This section is a quick recollection of the definition of virtual knot theory.
For more information the reader is referred to \cite{VKT,GPV,SVKT,DVK,Intro,KUP}.
The diagrammatic definition of virtual knot theory is that virtual knots and links are represented by diagrams that are like classical knot and link diagrams, except that there is added a new crossing called a {\it virtual crossing}. The virtual crossing is here indicated by a flat crossing (neither over nor under) that is encircled by a transparent circle. See Figure~\ref{virtualmoves} and Figure~\ref{detourmove}.
In these figures we indicate how the classical Reidemeister moves are generalized to include moves that  involve the virtual crossings. The general principle is that the virtual crossing behaves as an artifact of the projection of the virtual knot diagram to the plane. The actual virtual knot or link is independent of 
any embedding in the plane, but is assigned cyclic order of edges at each non-virtual crossing. Thus an attempt to embed the virtual knot in the plane can lead to extra crossings just as the embedding of a 
non-planar graph can require extra crossings. The moves for the virtual crossings are designed to respect this point of view. One can think of the virtual moves as generated by the local moves shown in 
Figure~\ref{virtualmoves} or one can say that they are generated by classical Reidemeister moves plus the detour move shown in Figure~\ref{detourmove}. The detour move allows an arc with a consecutive
sequence of virtual crossings to be excised and replaced any other such arc with consecutive virtual 
crossings.
\bigbreak

\begin{figure}
     \begin{center}
     \includegraphics[width=8cm]{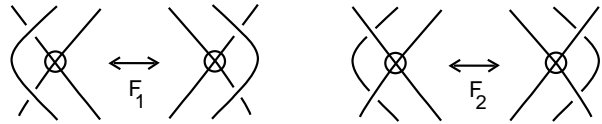}
     \caption{Forbidden Moves}
     \label{forbidden}
\end{center}
\end{figure}

The reader should take note of Figure~\ref{forbidden} where we illustrate two moves on virtual knots
that are not allowed and that do not follow from the given sets of moves. That these moves are forbidden is crucial in the structure of the theory. Allowing the move $F1$ we obtain {\it welded knot theory} a different variant that is closely related to the theory of welded braids of Rourke, Fenn and Rimiyani 
\cite{FRR}.
\bigbreak

Just as non-planar graphs may be embedded in surfaces of some genus, virtual knots and links can be represented by embeddings without virtual crossings in thickened orientable surfaces. In fact, the 
theory of virtual knots and links is equivalent to the theory of embeddings of circles in thickened surfaces modulo diffeomorphisms of the surfaces and one-handle stabilization of the surfaces. See
\cite{DVK,CS1,DK,KUP} for more information about this point of view.
\bigbreak

\subsection{Flat Virtual Knots and Links.}
{\it Flat virtual knots and links} are defined in exactly the same way as virtual knots and links except that the classical crossings are replaced by {\it flat crossings} indicated by transversely intersecting line segments and no information about over or under crossing. In the diagram of a flat virtual link one has 
flat classical crossings and the (also flat) virtual crossings. The moves are exactly the moves indicated in Figures~\ref{virtualmoves}, \ref{detourmove}, \ref{forbidden} with all crossings replaced by flat crossings. This means that once again we have forbidden moves, and it is easy to see that the presence of the forbidden moves gives flat virtual knot theory much non-trivial structure. In our construction of
the Affine Index Polynomial, we shall be first make an invariant labeling structure of flat virtual knots, and then use it to make the invariant for regular virtual knots.
\bigbreak

\section {The Polynomial Invariant}
\begin{figure}
     \begin{center}
     \includegraphics[width=10cm]{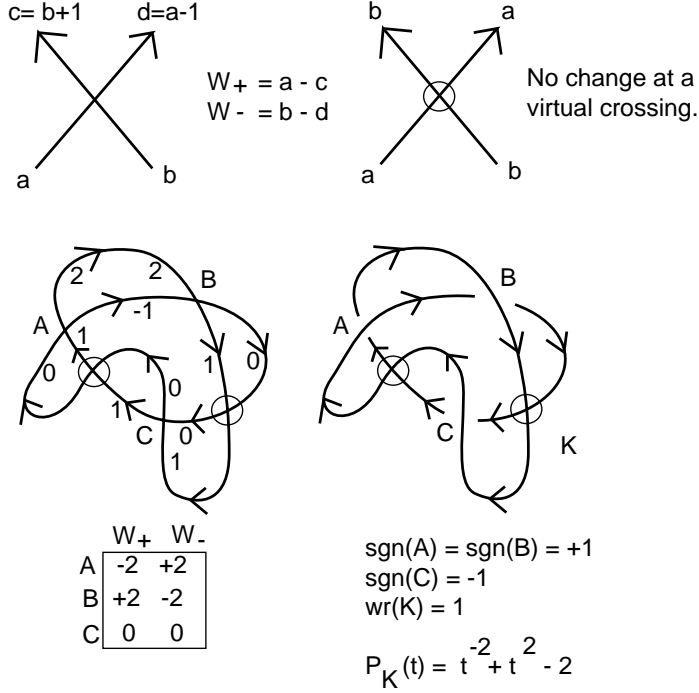}
     \caption{Labeled Flat Crossing and Example 1}
     \label{example1}
\end{center}
\end{figure}

We define a polynomial invariant of of virtual knots by first describing how to calculate the polynomial.
We then justify that this definition is invariant under virtual isotopy. Calculation begins with a flat oriented virtual knot diagram  (the classical crossings in a flat diagram do not have choices made for over or under). An {\it arc} of a flat diagram is an edge of the $4$-regualar graph that it represents. That is, an edge extends from one classical node to the next in orientation order. An arc may have many virtual crossings, but it begins at a classical node and ends at another classical node. We label each arc $c$ in the diagram with an integer $\lambda(c)$ so that an arc that meets  a classical node and crosses to the left increases the label by one, while an arc that meets a classical node and crosses to the right decreases the label by one. See Figure~\ref{example1} for an illustration of this rule. We will prove that such integer
labelling can always be done for any virtual or classical link diagram. In a virtual diagram the labeling is unchanged at a virtual crossing, as indicated in Figure~\ref{example1}. One can start by choosing some arc to have an arbitrary integer label, and then proceed along the diagram labelling all the arcs via this crossing rule.
We call such an integer labelling of a diagram a {\it Cheng coloring} of the diagram.
\bigbreak

\noindent {\bf Remark.} That any virtual knot diagram can receive a Cheng coloring is analogous to the fact that diagrams can be Fox-colored in $Z/NZ$ for appropriate values of $N.$ In the case of Fox colorings there is a corresponding algebra structure - {\it the quandle}. In the case of the Cheng coloring
we are looking at an example of a {\it flat affine biquandle}. We will discuss this algebraic background to the invariant in section 8 of this paper.
\bigbreak

\begin{figure}
     \begin{center}
     \includegraphics[width=6cm]{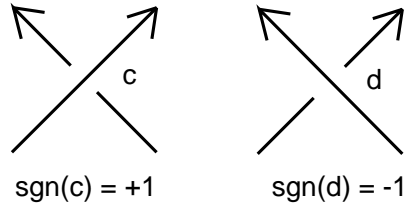}
     \caption{Crossing Signs}
     \label{crossingsign}
\end{center}
\end{figure}

Given a labeled flat diagram we define two numbers at each classical node $c$: $W_{-}(c)$ and 
$W_{+}(c)$ as shown in Figure 1. If we have a labeled classical node with left incoming arc $a$ and 
right incoming arc $b$ then the right outgoing arc is labeled $a - 1$ and the left outgoing arc is labeled
$b+1$ as shown in Figure~\ref{example1}. We then define
$$W_{+}(c) = a -( b +1)$$ and $$W_{-}(c) = b - (a -1)$$
Note that $$W_{-}(c) = - W_{+}(c)$$ in all cases.
\bigbreak

\noindent {\bf Definition.} Given a crossing $c$ in a diagram $K,$ we let $sgn(c)$ denote the sign of the crossing.
The sign of the crossing is plus or minus one according to the convention shown in 
Figure~\ref{crossingsign}.
The {\it writhe}, $wr(K),$  of the diagram $K$ is the sum of the signs of all its crossings.
For a virtual link diagram, labeled in the integers according to the scheme above, and a crossing 
$c$ in the diagram, define $W_{K}(c)$ by the equation
$$W_{K}(c) = W_{sgn(c)}(c)$$ where $W_{sgn(c)}(c)$ refers to the underlying flat diagram for $K$. Thus $W_{K}(c)$ is $W_{\pm}(c)$ according as the sign of the crossing is plus or minus. {\it We shall often indicate the weight of a crossing $c$ in a knot diagram $K$ by $W(c)$ rather than $W_{K}(c).$}
\bigbreak

Let $K$ be a virtual knot diagram. Define the {\it Affine Index Polynomial of $K$} by the equation
 $$P_{K} = \sum_{c} sgn(c)(t^{W_{K}(c)} - 1) = \sum_{c} sgn(c)t^{W_{K}(c)} - wr(K)$$ where the summation is over all classical crossings in the virtual knot diagram $K.$
We shall prove that the Laurent polynomial $P_{K}$
is a highly non-trivial invariant of virtual knots. 
\bigbreak

In Figure~\ref{example1} we show the computation of the weights for a given flat diagram and the computation of the polynomial for a virtual knot $K$ with this underlying diagram. The knot $K$ is an example of a virtual knot with unit Jones polynomial. The polynomial $P_{K}$ for this knot has the value
$$P_{K} = t^{-2} + t^{2} - 2,$$ showing that this knot is not isotopic to a classical knot. We will examine this example and others in the body of the paper.
\bigbreak

\section{Invariance of $P_{K}(t).$}
In order to show the invariance and well-definedness of  $P_{K}(t)$ we must first show the existence of labelings of flat virtual knot diagrams, as defined in section 1. We shall do this by showing that any virtual knot diagram $K$ that overlies a given flat diagram $D$ can be so labeled. 
\smallbreak

\noindent  {\bf Proposition.} Any flat virtual knot diagram has a Cheng coloring.
\smallbreak

\begin{figure}
     \begin{center}
     \includegraphics[width=8cm]{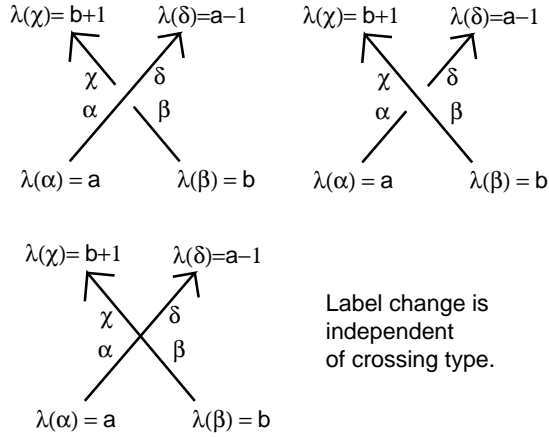}
     \caption{Labels for Crossings}
     \label{lambda}
\end{center}
\end{figure}

\noindent {\bf Proof.} Let $K$ be a flat virtual knot diagram.
Label the arcs $\alpha$ of $K$ by the following formula:
$$\lambda (\alpha) = \sum_{c \in O(\alpha)} sgn(c)$$
where $c$ denotes a classical crossing in $K,$ $sgn(c)$ is the sign of the crossing $c,$ and $O(\alpha)$ denotes the set of crossings first met as overcrossings on traversing, in the direction of its orientation, the diagram $K,$ starting at the arc $\alpha.$  The reader can easily check that this function 
$\lambda(\alpha)$ satisfies the formulas in Figure~\ref{example1} and Figure~\ref{lambda},  and hence can be computed for $K$ by finding the value for a single arc, and then using the information on the underlying flat diagram to obtain all the other values. The change of $\lambda(\alpha)$ across a classical crossing  does not depend upon the sign of the crossing, but only on whether one is crossing to the right or to the left. It is the fact that $\lambda(\alpha)$ is a signed 
sum of overcrossing encounters that makes its change independent of local crossing signs. This completes the proof of the Propposition. //
\bigbreak

Note also that if we follow the algorithm of the section 1, Figure~\ref{example1}, to compute a labeling, the resulting labeling will differ from the $\lambda$ labeling, defined above, by a constant integer at every label. Since the polynomial is defined in terms of the differences $W_{\pm}(c)$ at each classical crossing $c$ of $K,$ it follows that the weights $W_{\pm}$ as described in section 1, are well-defined.
\bigbreak

We can now prove a useful result about the weights.
Let $\bar{K}$ denote the diagram obtained by reversing the orientation of $K$ and let $K^{*}$ denote the diagram obtained by switching all the crossings of $K.$  $\bar{K}$ is called the {\it reverse} of $K,$ and $K^{*}$ is called the {\it mirror image} of $K.$
\smallbreak

\noindent {\bf Proposition.} Let $K$ be a virtual knot diagram and $W_{\pm}(c)$ the crossing weights
as defined in section 1 and above.  If $\alpha$  is an arc of $K,$ let $\bar{\alpha}$ denote the
corresponding arc of $\bar{K}$, the result of reversing the orientation of $K.$ 
\begin{enumerate}
\item Let $\lambda$ denote the label function, defined above, that counts 
overcrossings with signs. Then $$\lambda(\alpha) + \lambda(\bar{\alpha}) = wr(K)$$ where
$wr(K)$ denotes the writhe of the oriented diagram $K.$ 
\item Let $c$ be a crossing of $K$ and let $\bar{c}$
denote the corresponding crossing of $\bar{K}$, then $W(\bar{c}) = - W(c).$ 
\item Consequently, we have
$$P_{\bar{K}}(t) = P_{K}(t^{-1}).$$ Similarly, we have $$P_{K^{*}}(t) = -P_{K}(t^{-1}).$$
Thus this invariant changes $t$ to $t^{-1}$ when the orientation of the knot is reversed, and it 
changes global sign and $t$ to $t^{-1}$ when the knot is replaced by its mirror image.
\item If $K$ is a classical knot diagram, then for each crossing $c$ in $K$, 
$W(c) = 0$ and  $P_{K}(t) = 0.$ 
\end{enumerate}
\smallbreak

\noindent {\bf Proof.} The equation $\lambda(\alpha) + \lambda(\bar{\alpha}) = wr(K)$ follows immediately from the fact that crossing signs are not changed in reversal of orientation and that 
$O(\bar{\alpha}) = U(\alpha)$ where $U(\alpha)$ denotes the set of crossings first met as undercrossings on traversing, in the direction of the orientation of $\alpha.$ Thus
$$\lambda(\alpha) + \lambda(\bar{\alpha}) = O(\alpha) + U(\alpha) = wr(K)$$ since the
sum $O(\alpha) + U(\alpha)$ equals the writhe of $K$ for any arc in $K.$ From this fact and the definition of $W(c)$ as a difference of labels it is easy to calculate that $W(\bar{c}) = - W(c).$ The corresponding
equation $P_{\bar{K}}(t) = P_{K}(t^{-1})$ follows immediately from the definition of the polynomial.
In the case of the mirror image, the sign of each crossing changes, and weights change sign.
The formula $P_{K^{*}}(t) = -P_{K}(t^{-1})$ is a consequence of that. Finally, if $K$ is a classical knot diagram
and $k$ is a crossing in $K$ then $W(k)=0.$ The proof of this follows from the structure illustrated in Figure~\ref{classical}. In  that figure we remind the reader that the labeling of an arc adds one when the arc crossings another arc going to the right and subtracts one when the arc crosses another arc going to the left. We wish to prove that if we start by labeling a oriented arc at a crossing and then continue the labeling until we emerge at an outgoing arc from that crossing, then the outgoing arc will receive the same label as the ingoing arc. This will be true if and only if the number of right-going arcs crossed in the
loop are equal to the number of left-going arcs. First of all one can consider the loop that is traced in this journey. In a classical knot every crossing in that loop will be traversed twice, one with one handedness, and once with the opposite handedness. Thus there is zero total contribution from the self-crossings of the loop. What remains are the contributions from arcs that cross through the loop from the rest of the diagram. It is an easy Jordan curve argument to see that the total contributions from such curves is also zero (their self-crossings are irrelevant). The reader who examines Figure~\ref{classical} will be able to supply the remaining details that reduce the arguement to 
an application of the Jordan curve theorem. Since in a classical knot $K$ we now know that 
$W(k)=0$ for each crossing $k,$ it follows that $P_{K}(t)=0.$ This completes the proof. //
\bigbreak

\begin{figure}
     \begin{center}
     \includegraphics[width=8cm]{ClassicalWeights.EPSF}
     \caption{Classical Weights}
     \label{classical}
\end{center}
\end{figure}

\noindent {\bf Remark.} We will give examples where mirror images and reversals are distinguished by
$P_{K}(t)$ in the next section. The argument showing that $W(k)=0$ for crossings in a classical knot diagram can be used to understand when $W(k)$ is either zero or not zero in a classical crossing in a virtual diagram. We will point out relevant examples in the next section.
\bigbreak

\begin{figure}
     \begin{center}
     \includegraphics[width=8cm]{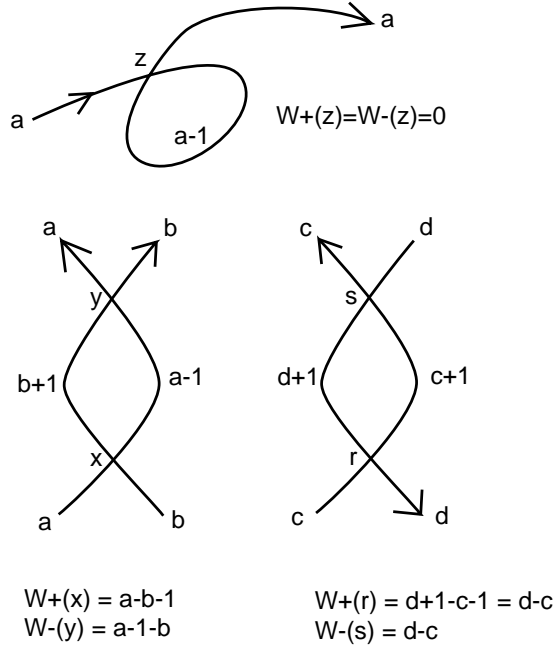}
     \caption{Reidemeister Moves 1 and 2}
     \label{r12}
\end{center}
\end{figure}

\begin{figure}
     \begin{center}
     \includegraphics[width=8cm]{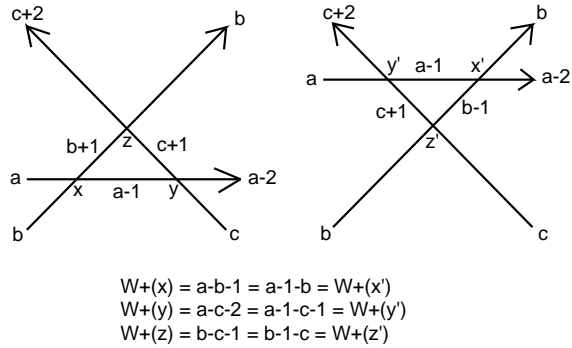}
     \caption{Reidemeister Move 3}
     \label{r3}
\end{center}
\end{figure}

We will now prove the invariance of $P_{K}(t)$ under virtual isotopy. The reader will recall that
virtual isotopy consists in the classical Reidemeister moves plus virtual moves that are all generated by one generic detour move. The (unoriented) virtual isotopy moves are illustrated in Figure~\ref{virtualmoves} and Figure~\ref{detourmove}. In Figure~\ref{r12} and Figure~\ref{r3} we show the relevant information for verifying that  $P_{K}(t)$ is an invariant of oriented virtual isotopy.
\bigbreak

\noindent {\bf Theorem.} Let $K$ be a virtual knot diagram. Then the polynomial $P_{K}(t)$
is invariant under oriented virtual isotopy and is hence an invariant of virtual knots.
\smallbreak

\noindent {\bf Proof.} Note that the definition of $P_{K}(t)$ makes it independent of the moves in
Figure~\ref{virtualmoves} involving virtual crossings.  The  labeling algorithm is independent of the purely virtual  moves and consequently the polynomial is invariant under them. Thus we need only verify invariance
under the standard oriented Reidemeister moves, shown as box $ $ in Figure~\ref{virtualmoves}.
Note that the Cheng coloring is uniquely inherited under the Reidemeister moves, and thus we only have to check the local changes in the coloring in relation to given types of move.
Since we wish to verify invariance under oriented Reidemeister moves, we use the well-known fact 
(see \cite{KP} page 81) that it is sufficient to verify invariance under type I moves, two orientations of type II moves and the single instance of the type III  move where there is a non-cyclic triangle in the center of the pattern and all the crossings have the same type (say positive). In Figure~\ref{r12} we give the relevant information for the moves of type I and II. For the type I move, we see that the weight $W_{\pm}(z)$ equals $0.$ This means that the contribution of a diagram with a type I move at a crossing $z$ is equal to $sgn(z).$ Since this is subtracted by the writhe $wr(K)$ in the formula $P_{K} = \sum_{c} sgn(c)t^{W_{K}(c)} - wr(K),$ we see that the polynomial is invariant under move I. For move II there are two cases or orientation as shown in the Figure~\ref{r12}. In each case, the move is available when the crossings have opposite sign. The calculation shown in the figure proves that when the crossings have opposite sign, the weights of the crossings are identical. Thus the sum of these two crossing contributions cancels out in the polynomial both in the $t$-terms and in the writhe term. Therefore the polynomial is invariant under the II move.
Finally, we examine the III move via the information in Figure~\ref{r3}. Here we assume that all three
crossings are positive. The labels are calculated from flat versions of the move. We see that
the resulting weights are simply permuted and the signs of the crossings remain unchanged. Thus the polynomial is invariant under move III. This completes the proof of the theorem. //
\bigbreak

\noindent {\bf Generalization from Knots to Links.} We are now in a position to generalize the invariant $P_{K}(t)$ to some cases of virtual links and even to some cases of classical links. It is possible that
a link diagram can be Cheng colored according to our rules. for example, view Figure~\ref{hlink} to see a labelling of the classical Hopf link. Before analyzing this figure, consider the proof we have given for 
the invariance of the polynomial $P_{K}(t).$ Cheng coloring is uniquely inherited under Reidemeister
moves and the weights at the three crossings of the third Reidemeister move are permuted under the 
move. These properties are true for the polynomial that we would write for any Cheng-colored link.
Thus we can conclude that if we are given pair $(L,C)$ where $L$ is a link diagram and $C$ is a 
Cheng-coloring of this diagram, then the polynomial $P_{L}(t)$, defined just as before, is an invariant
of the pair $(L,C)$ where a Reidemeister move applied to $(L,C)$ produces $(L',C')$ where $L'$ is the 
diagram obtained from $L$ by the move, and $C'$ is the coloring obtained from $C$ by the move.
\bigbreak

\section{Examples and Cobordisms}
In this section we will give a number of examples of computations of the Affine Index Polynomial.
\begin{enumerate}
\item We begin with the example in In Figure~\ref{example1}. Note since the Affine Index Polynomial
for this knot is $P_{K} = t^2 + t^{-2} -2,$ it follows that $K$ is not equivalent to its mirror image.
This knot $K$ is an example of a virtual knot with unit Jones polynomial. We refer the reader to examine 
\cite{VKT,Intro} for the details of this construction.
One can produce non-trivial virtual knots with unit Jones polynomial by the specific method of {\it virtualizing} a set of crossings in a classical knot  that would produce a classical unknot when switched.
See Figure~\ref{sw}. In this figure we illustrate {\it virtualization} of a crossing consisting in retaining it as an over or under crossing, but reversing its orientation. This is accomplished by redrawing the crossing, first scribing a virtual crossing, then performing the crossing with an opposite orientation, then scribing
another virtual crossing. The result is a crossing that is flanked by two virtual crossings such that smoothing the two virtuals gives the original diagram with a switched crossing. The figure also illustrates how the bracket polynomial of a virtualized crossing is the same as the bracket polynomial of the original diagram with a switched crossing. Given a non-trivial
classical knot $K$, one can choose a subset of crossings so that the unknot is obtained if they are all switched. If we virtualize this same set of crossings, we obtain a non-trivial \cite{VKT,Intro} virtual knot
$K'$ with unit Jones polynomial. {\it We conjecture that no such $K'$ is virtually isotopic to a classical 
knot.} This conjecture has so far only been verified in examples. The Affine Index Polynomial promises to be a useful tool in investigating this conjecture.
\bigbreak
 
 \begin{figure}
     \begin{center}
     \includegraphics[width=8cm]{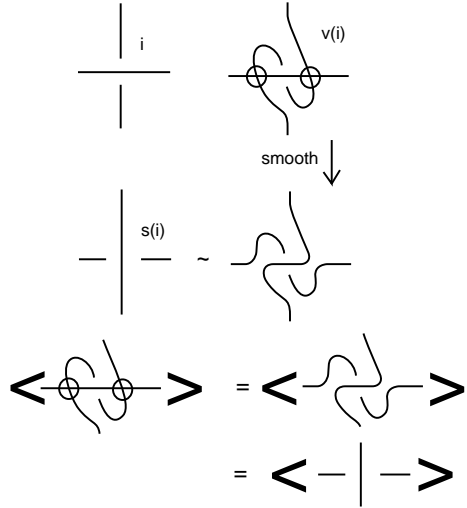}
     \caption{Switch and Virtualize}
     \label{sw}
\end{center}
\end{figure}

 \item Figure~\ref{family} Illustrates an infinite family of virtual knots with the same Affine Index Polynomial.
 \bigbreak
 
 \begin{figure}
     \begin{center}
     \includegraphics[width=8cm]{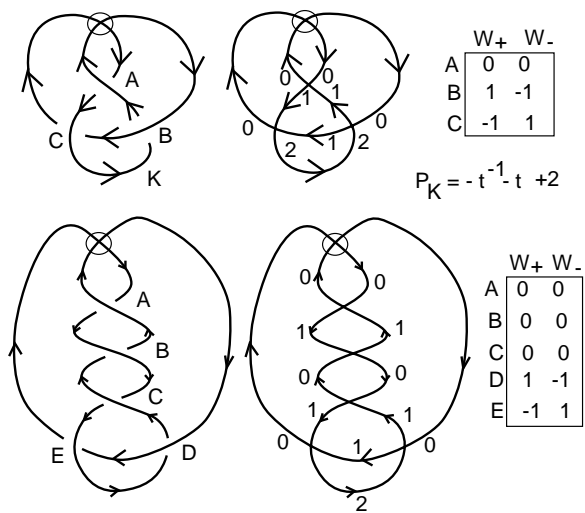}
     \caption{A Family of Virtual Knots with the Same Polynomial}
     \label{family}
\end{center}
\end{figure}

 \item Figure~\ref{virtex} gives another example of virtualization. The knot $K$ has unit Jones polynomial
 but a non-trivial Affine Index Polynomial, proving that $K$ is not classcial. Then in Figure~\ref{basicob} and 
 Figure~\ref{cob} we illustrate how the appearance of zeroes in the list of vertex weights for the polynomial can be used to produce labelled knots and links where the crossings with null weights have been smoothed.  We will call the smoothing indicated in  Figure~\ref{basicob} a {\it basic labeled cobordism}. Thus if a knot has crossings with null weights, then it is labeled cobordant to a link with only non-zero weights (or an empty set of weights). While not all links can be labeled, this form of cobordism does produce labeled links, and the Index Invariant can be extended to such links as indicated in 
 Figure~\ref{hlink}. Here we write down the most general labeling for the link, and then deduce a set of variable integer exponents for the polynomial invariant. We shall leave the details of this generalization 
 to another paper.
 \bigbreak
 
 \begin{figure}
     \begin{center}
     \includegraphics[width=8cm]{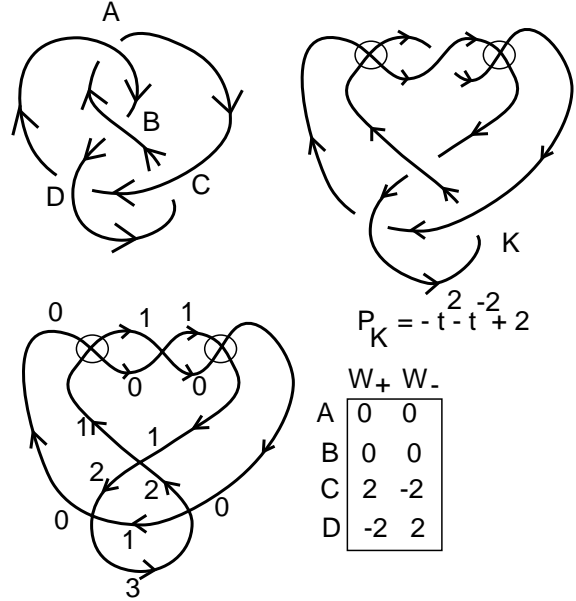}
     \caption{A Virtualization Example}
     \label{virtex}
\end{center}
\end{figure}

\begin{figure}
     \begin{center}
     \includegraphics[width=8cm]{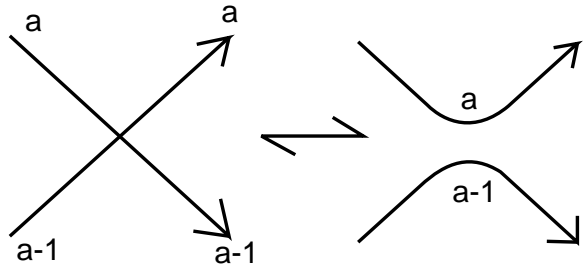}
     \caption{Basic Labeled Cobordism}
     \label{basicob}
\end{center}
\end{figure}

\begin{figure}
     \begin{center}
     \includegraphics[width=8cm]{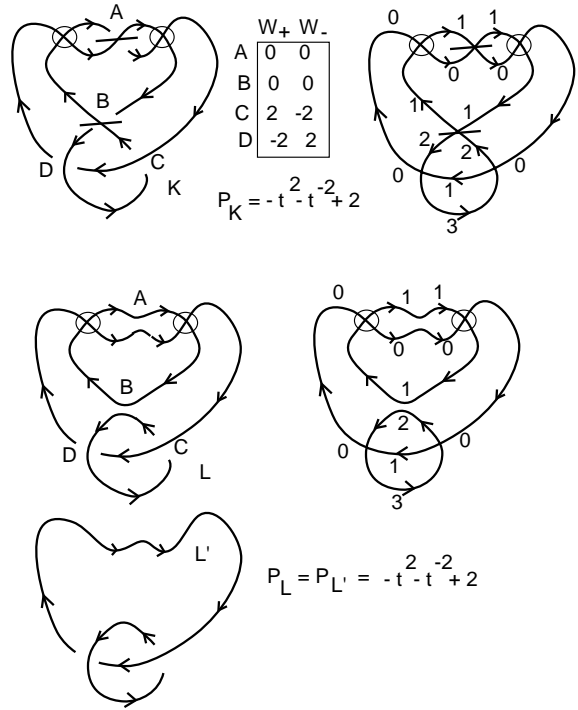}
     \caption{Labeled Cobordism of a Knot to a Link}
     \label{cob}
\end{center}
\end{figure}

\begin{figure}
     \begin{center}
     \includegraphics[width=8cm]{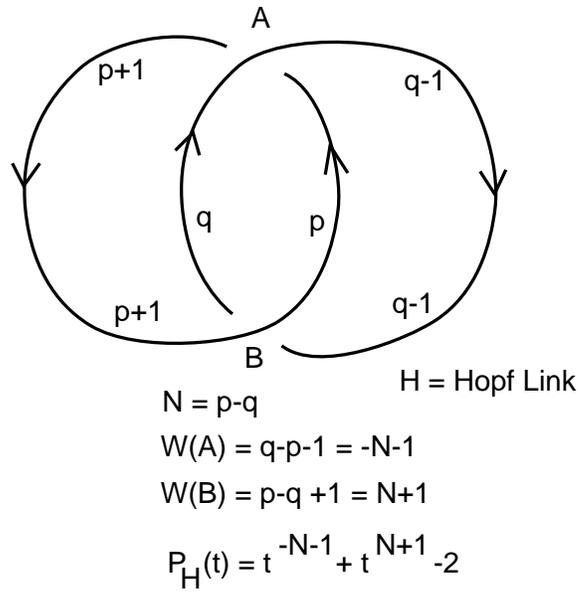}
     \caption{Invariant for the Hopf Link}
     \label{hlink}
\end{center}
\end{figure}

 \item In Figure~\ref{threeswitch} we give an example of a classical knot that can be transformed to an 
 unknot by switching three crossings. This knot has unknotting number three. However when we take the corresponding virtualization, we find the the Index invariant is equal to zero. Thus the non-classicality of this virtualization is not detected by the Index invariant. This leads to a question that is worth investigating: {\it Characterize those virtualizations that are detected by the Index invariant.}
 \bigbreak
 
 \begin{figure}
     \begin{center}
     \includegraphics[width=8cm]{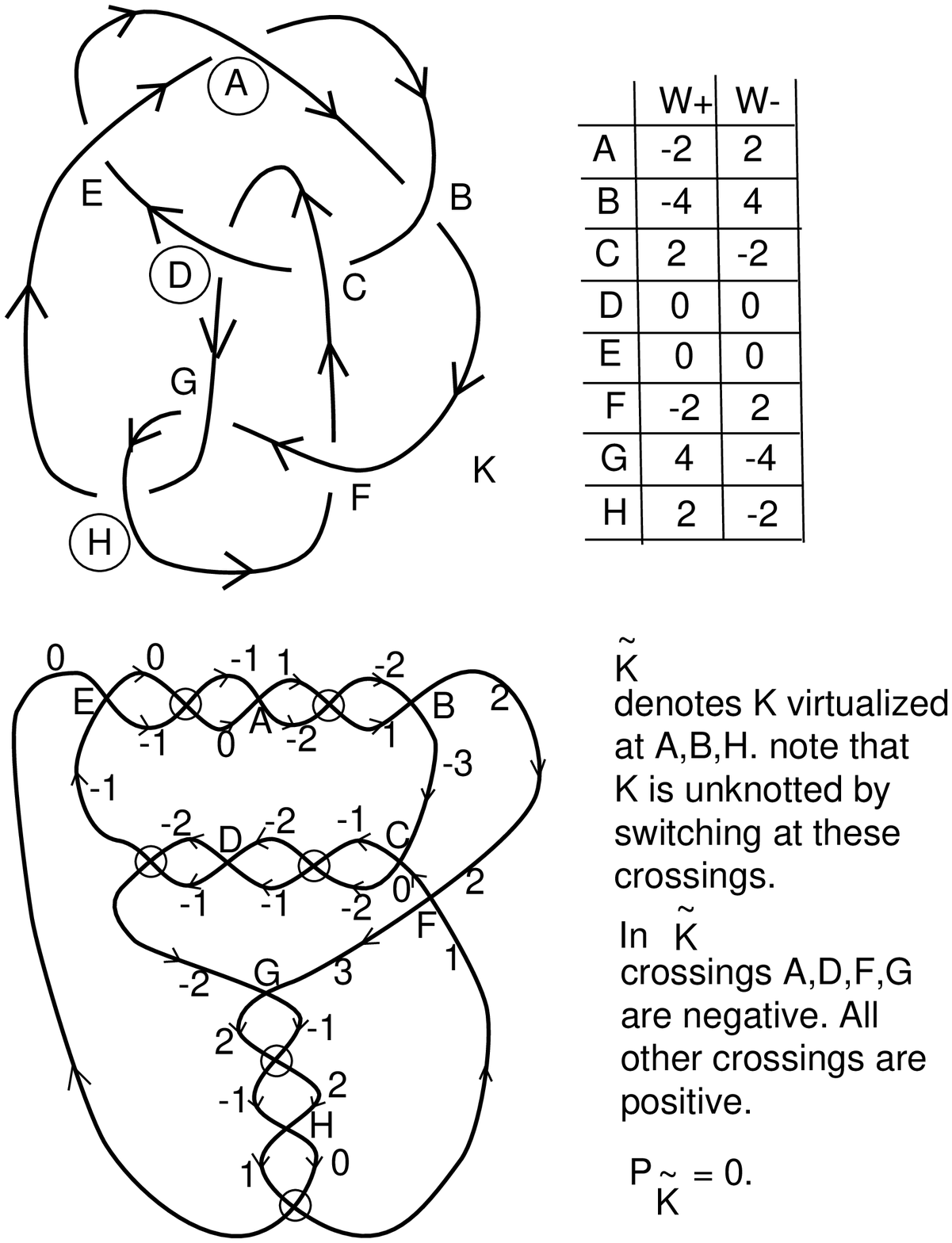}
     \caption{A Virtualization with Trivial Polynomial}
     \label{threeswitch}
\end{center}
\end{figure}

 \item In Figure~\ref{flatnontriv} we give an example of a flat diagram such that the Index invariant is not
 zero for any choice of resolution for its crossings. This implies that the flat diagram $D$ is itself a non-trival virtual flat becuase it is not hard to see that if there were a flat isotopy that trivializes $D$ then it would be overlaid by an trivializing isotopy for some choice of crossings for the flat diagram. This example shows that one can sometimes use the Index invariant to detect non-trivial flat knots.
 In this last example we see from the Index invariant that all knots overlying this flat diagram are 
 non-classical, non-invertible and inequivalent to their mirror images.
 \bigbreak

\begin{figure}
     \begin{center}
     \includegraphics[width=8cm]{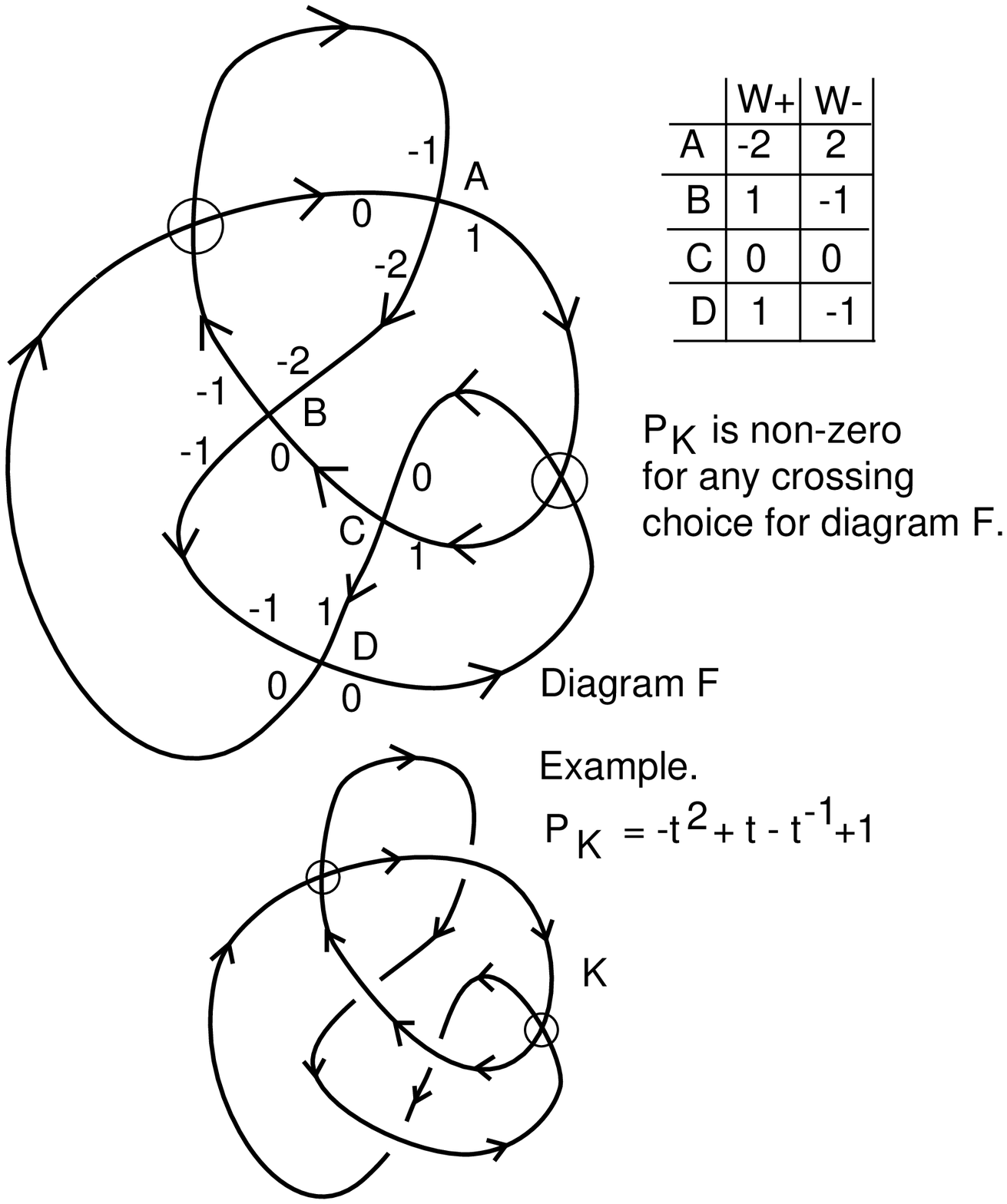}
     \caption{A Non-trivial and Non-invertible Flat Knot}
     \label{flatnontriv}
\end{center}
\end{figure}

 \item In Figure~\ref{zero} we illustrate a non-trivial knot $KZ$ with trivial Affine Index Polynomial. All the crossing weights are zero for this knot.  Its non-triviality can be checked by computing its Jones polynomial. But the knot is not shown to be non-classical by the Arrow Polynomial \cite{DKArrow},
Arrow Polynomial Categorification \cite{DyeKM,KaestKauff},  or by the Sawollek Polynomial \cite{KR}. At this  writing we do not know whether $KZ$ is non-classical but we conjecture that this is the case.  It is not hard to make infinitely many examples of this kind (for example, by adding more twists to the given example) and so we are led to search for new invariants to detect
 non-classicality. A study of examples of this kind will be the subject of another paper.

 \noindent {\bf Note added in proof.} The example in Figure~\ref{zero} has been shown to be non-classcial by M. V.  Zenkina, using results in her paper \cite{Zenkina} and it has been shown to be non-classcial by the Author and Slavik Jablan by using results in the paper by V. O.  Manturov \cite{ParityProject}. In the latter case, the relevant fact is Manturov's Theorem that states the no classical knot has virtual representatives with the number of classical crossings less than its minimal crossing number as a classical knot. 

\begin{figure}
     \begin{center}
     \includegraphics[width=8cm]{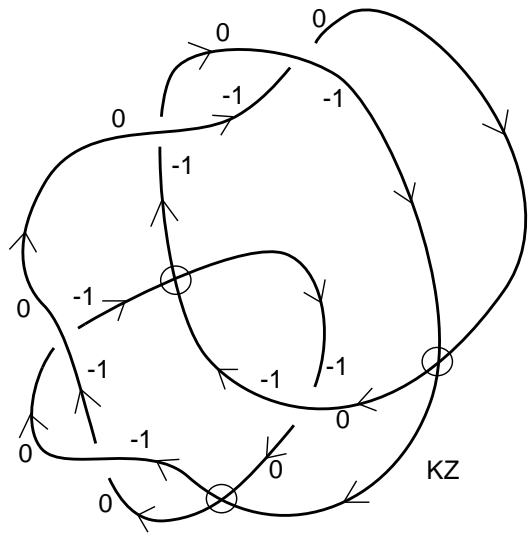}
     \caption{A Non-Trivial Knot With Zero Affine Index Polynomial}
     \label{zero}
\end{center}
\end{figure}

\end{enumerate}

\section{Vassiliev Invariants From the Affine Index Polynomial}
The method of exponential substitution yields an infinite set of Vassiliev invariants
for the Affine Index Polynomial. We show how to do this in the present section, and give specific formulas for these invariants in terms of the weights $W_{K}(c).$
\bigbreak

Note that we have the following skein relation for the Affine Index Polynomial:
$$P_{K+} - P_{K-} = t^{W_{+}(c)} + t^{W_{-}(c)} -2 =  t^{W_{+}(c)} + t^{-W_{+}(c)} -2$$
where $c$ denotes the crossing at which the switch occurs, \
$W_{+}(c) = W_{K{+}}(c), W_{-}(c) = W_{K{-}}(c)$ and $K_{+}$ is a diagram with a positive crossing at $c$ and $K_{-}$ is the diagram with a negative crossing at $c$ that is obtained by 
changing just this one crossing. We leave the proof of this identity to the reader. Note that the $-2$ comes from the writhe difference between the two diagrams.
\bigbreak

We use the skein difference to extend the Affine Index Polynomial to an invariant of virtual $4$-regular graph embeddings (See \cite{VKT}) by defining $$P_{K_{*}} = P_{K+} - P_{K-} $$ where $K_{*}$ denotes the
original diagram with a graphical node at the site $c.$ This formula then defines the graph invariant
(by expansion into differences) for any virtual graph $G.$
Let $t = e^{x}.$ Then $$P_{K}(e^{x}) = v_{0}(K) + v_{1}(K)x + v_{2}(K)x^{2} + \cdots$$
where $v_{i}(K)$ is the coefficient of $x^{i}$ in this power series. Note also, from the skein relation above, that $$P_{K_{*}}(e^{x})  =  (e^{xW_{+}(c)} + e^{-xW_{+}(c)} - 2).$$
Since
$$x^{2} | (e^{xW} + e^{-xW} - 2)$$ for any integer $W$, this implies that 
$x^{2n} | P_{G}$ whenever $G$ has more than $n$ nodes. This means that 
$v_{k}(G)=0$ whenever $n \ge k/2.$ In turn, this means that $v_{k}(G)$ is a Vassiliev
invariant of order $\lceil k/2 \rceil.$ (A Vassiliev invariant is of order $m$ if it vanishes on graphs with 
more than $m$ nodes.)
\bigbreak

To obtain formulas for these Vassiliev invariants, consider any Laurent polynomial
$$P(t) = a_{1}t^{c_{1}} + \cdots + a_{m}t^{c_{m}}$$
where the exponents are an increasing set of integers $\{ c_{1} < \cdots <  c_{m} \}$ and the coefficients
$a_{k}$ are also integers. Then 
$$P(e^{x}) = \sum_{i=1}^{m} a_{i}e^{c_{i}x} = \sum_{i=1}^{m} a_{i} \sum_{k=1}^{\infty} c_{i}^{k} x^{k}/k!
=  \sum_{k=1}^{\infty} (1/k!)(\sum_{i=1}^{m} a_{i}c_{i}^{k}) x^{k}.$$
Thus the coefficient of $x^{k}$ in the exponential substution is 
$$v_{k} = (1/k!)(\sum_{i=1}^{m} a_{i}c_{i}^{k}).$$  
We have proved
\smallbreak

\noindent {\bf Proposition.} For the Index invariant, we have a finite type Vassiliev invariant 
of order $\lceil n/2 \rceil$ given by the formula
$$v_{n}(K) = (1/n!)(\sum_{c} sgn(c) W_{K}(c)^{n})$$  where $c$ runs over all the crossings of the knot 
$K$, $sgn(c)$ is the sign of the crossing, and $W_{K}(c)$ is the index of the crossing.
\bigbreak

\noindent {\bf Remark.}
\begin{enumerate}
\item While we have used the polynomial to deduce these Vassiliev invariants, it is worth noting that 
they depend only on the signs and weights of the crossings.
\item Note that since $W_{\bar{K}}(c) = - W_{K}(c)$ it follows that if $v_{n}(K)$ is non-zero for any
odd $n,$ then $K$ is inequivalent to its reverse orientation $\bar{K}.$
\item It is easy to see that $v_{1}(K) = 0$ for any virtual knot $K,$ since the sum of all the signed weghts
$sgn(c)W_{K}(c)$ is equal to the sum of all the $W_{+}(c)$ across the underlying flat diagram. And this latter sum is zero because every edge label occurs with both a positive and a negative sign in the sum.
However, $v_{3}(K)$ (a Vassiliev invariant of order $2$) is often non-zero in many examples, showing that irreversibility can be detected at the level of a Vassiliev invariant of order two.
\item A specific example that has this property is shown in Figure~\ref{flatnontriv} with the specific knot $K$ with polynomial $P_{K} = -t^2 + t - t^{-1} + 1.$ We then have the formula
$$v_{n}(K) = (1/n!)(- 2^n + 1^n - (-1)^n) = (1/n!)(1 + (-1)^{n+1} - 2^n ) .$$
Thus $v_{1}(K) = 0$ but $v_{3}(K) = (1/3!)(2 - 2^3) = -1,$ proving (again) that $K$ is inequivalent to its
reverse.
\end{enumerate}
\bigbreak

\section{Flat Biquandles}
We have seen that the Affine Index Polynomial is based upon a labeling of a flat virtual diagram by integers that follows the rule given in Figure~\ref{lambda}. The form of this labeling, following the rules of this 
figure, is invariant under virtual isotopy in the sense that for each virtual move, there is a unique way to extend the labelling to the new diagram, changing it only locally at the site of the move. It is the purpose of this section to give the Affine Index Polynomial a context by discussing the algebraic structure behind this
labeling. To this purpose we make the following definition.
\bigbreak

\begin{figure}
     \begin{center}
     \includegraphics[width=8cm]{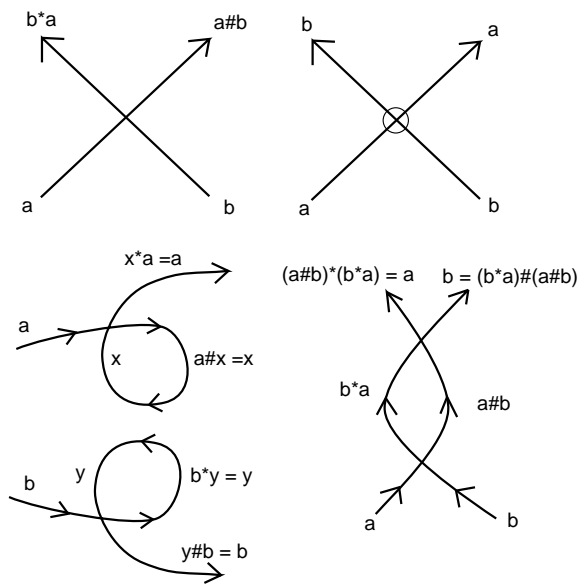}
     \caption{Flat Biquandle Operations and Moves 1 and 2.}
     \label{flat1}
\end{center}
\end{figure}

\begin{figure}
     \begin{center}
     \includegraphics[width=6cm]{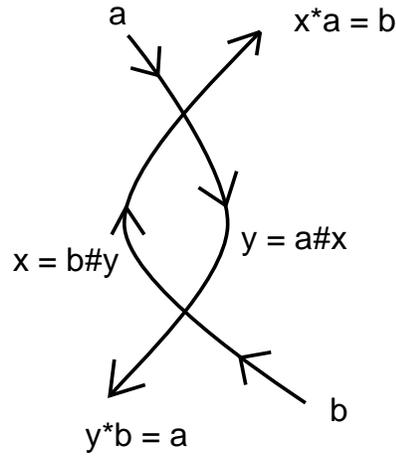}
     \caption{Flat Biquandle Reverse Move 2.}
     \label{flat2}
\end{center}
\end{figure}

\begin{figure}
     \begin{center}
     \includegraphics[width=8cm]{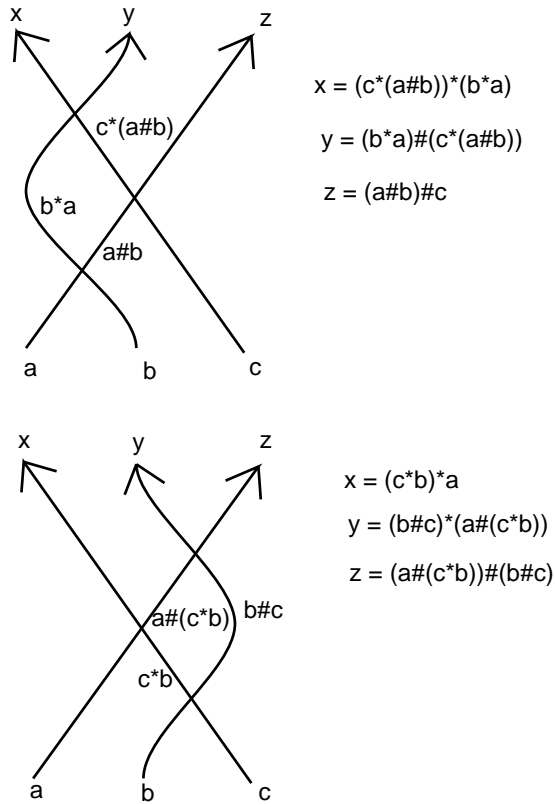}
     \caption{Flat Biquandle Move 3.}
     \label{flat3}
\end{center}
\end{figure}

\noindent {\bf Definition.} A {\it flat biquandle} is an algebraic system $S$ with two binary operations
denote $a*b$ and $a\sharp b$ satisfying the following axioms.
\begin{enumerate}
\item For each $a \in S$ there is a unique $x \in S$ such that $a\sharp x = x$ and $x*a = a.$
And for each  $a \in S$ there is a unique $y \in S$ such that $a*y = y$ and $y \sharp a = a.$
\item For all $a,b \in S,$ $$(a \sharp b)*(b*a) = a$$ and $$(b*a)\sharp (a \sharp b) = b.$$
For all $a,b \in S$ there exist unique elements $x,y \in S$ such that $$x = b \sharp y,$$
$$y = a \sharp x,$$ $$b = x*a,$$  $$a = y*b .$$
\item For all $a,b,c \in S$ 
$$(a \sharp b) \sharp c = (a \sharp (c*b)) \sharp (b \sharp c),$$
$$(c*b)*a = (c*(a \sharp b))*(b*a),$$
$$(b \sharp c)*(a \sharp (c*b)) = (b*a) \sharp (c*(a \sharp b)).$$
\end{enumerate}
\bigbreak 

The motivation for this definition is supplied in Figure~\ref{flat1}, Figure~\ref{flat2} and Figure~\ref{flat3}.
In Figure~\ref{flat1} we interpret the operations $a*b$ and $a \sharp b$ in terms of labeling a flat diagram at a classical flat crossing. At a virtual crossing, there is no change in the labeling for arcs going
through the virtual crossing. In this figure we show how the relation $1$ above corresponds to the 
flat first Reidemeister move and we show how the first part of relation $2$ above corresponds to the
oriented flat second Reidemeister move. In Figure~\ref{flat2} we show how the second part of 
relation $2$ corresponds to the reverse oriented second Reidemeister move. Finally, in 
Figure~\ref{flat3} we show how relation $3$ above corresponds to invariance of labeling under the 
flat third Reidemeister move. This means that if a flat virtual diagram is consistently colored with the elements of a flat biquandle, then this coloring can be uniquely modified to produce colorings on diagrams that are equivalent to the given diagram under flat virtual equivalence. And also, we can define the flat biquandle $FB(K)$ for a virtual diagram $K$ by taking the free flat biquandle (in the sense of universal algebra) on labels for the arcs of the diagram (extending from classical crossing to classical crossing) modulo the relations introduced by the definition of operations in Figure~\ref{flat1}.
Flat biquandles have been studied in \cite{HenrichNelson} and are called there {\it semiquandles}. 
We shall call an algebraic system that satisfies only $1.$ and $2.$ above a {\it flat pre-biquandle}.
 and for short we will use the term {\it preflat}. If $S$ is a preflat then it can be used to examine invariance
 under the first and second flat Reidemeister moves.
\bigbreak

It is our purpose here to introduce flat biquandles  in  relation to the 
Affine Index Polynomial invariant for virtual knots and links. In particular, it is easy to see that the following structure on the integers $Z$ is a flat biquandle: $$a*b =  a + 1,$$ $$a \sharp b = a - 1.$$ This is exactly the underlying label rule for the Index polynomial. Thus we see that the underlying structure for this polynomial invariant is a very specific flat biquandle structure on the integers, that allows labeling of 
any flat virtual knot diagram. In this light of this, we wish to determine when a flat biquandle structure can be given by an {\it affine formula} of the following type:
$$a*b = ra + sb + k,$$
$$a \sharp b = pa + qb + l.$$
Here we assume that the underlying set $S$ of the flat biquandle is a module over a commutative ring $R$ with $r,s,p,q$ elements of $R$ and $k,l$ specific elements of $S.$ If these formulas define the structure of a flat biquandle on the set $S,$ then we say that $S$ is an {\it affine biquandle}.
\bigbreak

\noindent {\bf Theorem.} Let $S$ be a module over a commutative ring $R$ with unit $1 = 1_{R}$ and
no zero divisors.  Let $k$ be any element of $S$ and let $\alpha$  an invertible element of $R.$  Then the following formulas define the general affine linear flat biquandle with coefficients in $R$
$$a*b  =  a^{*}  = \alpha a + k$$ and
$$a \sharp b = a^{\sharp} =  \alpha^{-1} a - \alpha^{-1} k$$
When $\alpha =1$ this flat biquandle is the root structure of the Index polynomial studied in the present paper. Note that we have chosen the notations $a^{*} ,a^{\sharp}$ to indicate that the operations in 
this flat biquandle are actually unary operations on the set $S.$
\bigbreak

\noindent {\bf Proof.}  We begin with 
$$a*b = ra + sb + k$$ and
$$a \sharp b = pa + qb + l.$$
We first look at the pair of equations
$$a \sharp x = x$$ and 
$$x*a = a.$$ We then have
$$ x = a \sharp x = pa + qx + l$$ and
$$a = x*a = rx + sa + k.$$
From the first equation we have
$$(1-q)x = pa + l.$$ From this we conclude that we shall need that $1-q$ is invertible in $R$ and 
that 
$$x = \frac{p}{1-q}a + \frac{l}{1-q}.$$
Putting this into the second equation, we have
$$a = \frac{rp}{1-q}a + \frac{rl}{1-q} + sa + k.$$
Thus
$$(\frac{rp}{1-q} + s -1)a + (\frac{rl}{1-q} + k) = 0.$$
Since this must be true for all $a$, we require that the two coefficients vanish. Rewriting, we 
have
$$rp + (1-q)s -(1-q) = 0$$ and
$$rl + k(1-q) = 0.$$
\bigbreak

Now we turn to the second flat Reidemeister move.This corresponds to the equations
$$a = (a \sharp b)*(b*a)$$
and
$$b = (b*a)\sharp (a\sharp b).$$
A calculation shows that
$$(b*a)\sharp (a\sharp b) = (pr +q^2 )b + p(s + q)a + (pk + ql + l).$$
Thus we demand that 
$$pr + q^2 = 1$$
$$p(s + q) = 0$$
$$pk + ql + l = 0.$$
The analogous calculation with the other equation from the second move yields the further 
conditions
$$pr + s^2 = 1$$
$$r(s + q) = 0$$
$$rl + sk + k = 0.$$
Using our assumption that the ring $R$ does not have zero-divisors, we conclude that
either $q = -s$ or $p = r = 0.$ Consider the case where $p = r = 0.$  It is easy to see that in this case the resulting operations are not invariant under the reversed orientation form of the second Reidemeister move. We leave this
for the reader to verify. Therefore, we shall assume that $s = -q$ and find that 
the equations reduce to 
$$rp = 1 - q^2 $$
$$pk + (q +1) l = 0$$
$$rl +k(1-q) = 0.$$
with 
$$a*b = ra -qb + k$$
$$a \sharp b = pa + qb + l.$$
\bigbreak

Now we apply the further condition of invariance under the reverse oriented Reidemeister two move.
This is the requirement that there exist unique elements $x,y \in S$ such that $x = b \sharp y$,
$y = a \sharp x$, $x * a = b$ and $a= y * b.$ Thus we have
$$a = y * b = ry -qb + k$$
$$b = x * a = rx - qa + k.$$
Hence
$$y = \frac{qb + a - k}{r}$$
$$x = \frac{qa + b -k}{r}.$$
Using these forms for $x$ and $y$, we have
$y = a \sharp x$ as the equation
$$ pa + q\frac{qa + b -k}{r} + l = \frac{qb + a -k}{r}.$$
This is equivalent to the condition
$$(rp - 1 + q^2 )a -qk + rl + k = 0.$$
Similarly, $x = b \sharp y$ yields the condition
$$(rp -1 +  q^2 )b -qk + rl +k = 0.$$
Thus we conclude that the following conditions will ensure invariance under the first and both of the second flat Reidemeister moves:
$$rp + q^2 = 1$$
$$pk + l(q +1) = 0$$
With $p \ne 0$ and $1+q$ invertible, we have
$$r = \frac{1-q^2}{p}$$
$$l =\frac{-p}{1 +qk}$$
so that 
$$a*b =(\frac{1-q^2}{p})a -qb + k$$ and
$$a \sharp b = pa + qb - (\frac{p}{1 +q})k.$$
This is our final set of general equations for a possible linear affine flat biquandle.
It then remains to check invariance under the flat third Reidemeister move. 
We checked this invariance by using a small computer program and found that it demands that
$q =0.$ Thus we are left with the equations
$$a*b = p^{-1}a + k$$ and
$$a \sharp b = pa - pk.$$
This completes the proof of the theorem.//
\bigbreak

\begin{figure}
     \begin{center}
     \includegraphics[width=6cm]{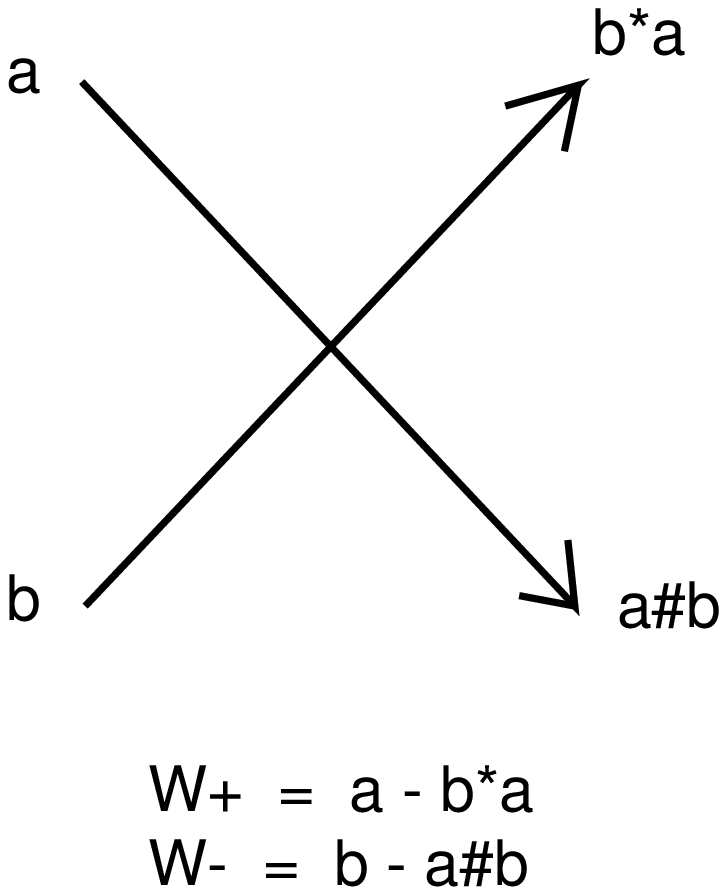}
     \caption{Weights from a Flat Biquandle.}
     \label{flat4}
\end{center}
\end{figure}

\noindent{\bf Remark.} Our purpose in this section  is put the index polynomial into the context of the flat biquandle.  To this end, lets determine when one can use this flat affine biquandle of the above theorem to produce a weight system for a polynomial invariant.  In Figure~\ref{flat4} we have illustrated a crossing with biquandle labels. We see from the figure that generalized weights at this crossing would be
$$W_{+} =  a - b*a$$ and 
$$W_{-} =  b - a \sharp b.$$
In order to generalize our results that gave an invariant polynomial from the weight system, we need that
$$W_{+} + W_{-} = 0.$$
\bigbreak

\noindent {\bf Lemma.} In order to have $W_{+} + W_{-} = 0,$  with in a lableling by a general preflat, we need that $$a + b =  b*a +  a \sharp b.$$
\smallbreak

\noindent {\bf Proof.} This follows immediately from the discussion above. //
\bigbreak

\noindent {\bf Lemma.} Let the preflat defined by the equations below be called the
{\it basic affine preflat}:
$$a * b = (1-q)a - qb + k$$
$$a \sharp b = (1 + q)a + qb - k.$$
This structure is a preflat by the argument in the last theorem (let $p = 1+q$),
and the weights defined as
$$W_{+} =  a - b*a$$ and 
$$W_{-} =  b - a \sharp b.$$
satisfy the condition 
$$W_{+} + W_{-} = 0.$$ Thus the basic affine preflat can be used to define invariants of flat oriented virtual diagrams up to the equivalence relation generated by the first and second flat Reidemeister moves.
\smallbreak

\noindent {\bf Proof.}
We find
$$b*a +  a \sharp b = (1-q)b - qa + k + (1+q)a + qb -k$$
$$= a + b.$$
Thus we conclude that this preflat will yield a weight system for any value of $q.$
This completes the proof of the Lemma. //
\bigbreak

\
We leave it to the reader to check that for the more general preflat in the proof of the theorem, the restriction added by requiring that $$W_{+} + W_{-} = 0$$ gives exactly that
 $p = 1+q$ and $r = 1-q.$)
 It follows that the basic affine preflat can be used to create pre-invariants of virtual knots where a 
 {\it pre-invariant} is invariant under the first and second Reidemeister moves, but not invariant under the third flat Reidemeister move. The pre-invariant can be written formally as
$$P_{K} = \sum_{c} sgn(c) t^W_{K}(c) - wr(K)$$ 
where $t$ is a formal variable  and $W_{K}(c)$ denotes the weight of the crossing defined as above for a 
labeling of the knot by the preflat. This invariant also depends upon the choice of coloring of the 
knot diagram by elements of the preflat. To make an invariant that is dependent only on the diagram, one can form a further sum over all such colorings. There is more work to be done in this domain.
\bigbreak

\noindent {\bf Remark.} We can call flat virtual diagrams taken up to all oriented versions of the flat first and second Reidemeister moves {\it virtual doodles} following \cite{FennDoodle,KhoDoodle}. In these 
papers by Fenn and Khovanov one considers flat link diagrams (not virtual) taken up to the equivalence relation generated by the first two Reidemeister moves. Study of virtual doodles will be undertaken in a separate paper.
\bigbreak

\noindent {\bf Remark.} In the theorem, we arrived at the affine biquandle given by the unary operations
$$ a*b  = a^{*} = \alpha a + k $$ and
$$ a \sharp b = a^{\sharp} =  \alpha^{-1} a - \alpha^{-1} k .$$
It is easy to see that if a flat biquandle is defined by unary operations
$$a*b  = a^{*}$$  and
$$a \sharp b = a^{\sharp},$$
then it suffices that these operations be inverses of one another with  
$$(a^{\sharp})^{*} = (a^{*})^{\sharp} = a$$
for the structure to define a flat biquandle. We leave this verification to the reader.
\bigbreak

Note that for this simple affine flat biquandle we have
$$W_{+} + W_{-}  =  -b*a  - a\sharp b +a +b =  -b^{*} - a^{\sharp} + a + b.$$ 
Thus for the operations above, we have
$$W_{+} + W_{-}  =  -\alpha b - k - \alpha^{-1} a +\alpha^{-1} k + a + b.$$ 
We conclude that $W_{+} + W_{-}  = 0$ exactly when $\alpha = 1.$ This shows that for a direct 
invariant there is no choice other then the very simple flat biquandle structure on which our Index polynomial is based.
\bigbreak



\end{document}